\documentclass[12pt,reqno]{article}
\usepackage{fullpage}
\usepackage{enumerate}
\usepackage{upgreek}
\usepackage{booktabs}
\usepackage{amsthm,amsmath,amssymb}
\usepackage{tikz}
\usepackage{subfigure}
\usepackage{graphicx}
\usepackage[export]{adjustbox}
\definecolor{webgreen}{rgb}{0,.5,0}
\definecolor{webbrown}{rgb}{.6,0,0}
\usepackage[colorlinks=true,
linkcolor=webgreen,
filecolor=webbrown,
citecolor=webgreen]{hyperref}

\newcommand{\seqnum}[1]{\href{http://oeis.org/#1}{\underline{#1}}}
\newcommand{\Figs}[1]{\hyperref[#1]{Figure~\ref*{#1}}}
\newcommand{\Tabs}[1]{\hyperref[#1]{Table~\ref*{#1}}}
\newcommand{\Lem}[1]{\hyperref[#1]{Lemma~\ref*{#1}}}
\everymath{\displaystyle}

\theoremstyle{plain}
\newtheorem{theorem}{Theorem}
\newtheorem{lemma}[theorem]{Lemma}
\newtheorem{proposition}[theorem]{Proposition}

\theoremstyle{remark}
\newtheorem{remark}[theorem]{Remark}

\title{\bf  A generating polynomial for the two-bridge knot with Conway's notation $ \mathcal{C}(n,r) $}
\author{Franck Ramaharo\\
\small D\'epartement de Math\'ematiques et Informatique\\[-0.8ex]
\small Universit\'e d'Antananarivo\\[-0.8ex] 
\small 101 Antananarivo, Madagascar\\
\small\href{mailto:franck.ramaharo@gmail.com}{\tt franck.ramaharo@gmail.com}\\
}

\date{\small\today\\}

\begin{document}
\maketitle		

\begin{abstract}
We construct an integer polynomial whose coefficients enumerate the Kauffman states of  the two-bridge knot  with Conway's notation $ \mathcal{C}(n,r) $. 

\end{abstract}

\section{Introduction}
\begin{figure}[!ht]
\centering
\includegraphics[width=.31\linewidth]{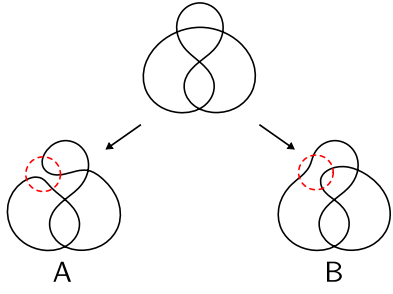}
\caption{The states of a crossing.}
\label{Fig:figure8}
\end{figure}

A \textit{state} of a knot shadow diagram is a choice of splitting its crossings \cite[Section 1]{Kauffman1}. There are two ways of splitting  a crossing:
\[
(\mathsf{A})\  \protect\includegraphics[angle=90,origin=c,width=0.06\linewidth,valign=c]{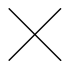}\Longrightarrow\protect\includegraphics[angle=90,origin=c,width=0.06\linewidth,valign=c]{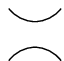},
\qquad(\mathsf{B})\
\protect\includegraphics[width=0.06\linewidth,valign=c]{crossing}\Longrightarrow\protect\includegraphics[width=0.06\linewidth,valign=c]{A_split}.
\]
By \textit{state of a crossing} we understand either of the split of type $ (\mathsf{A}) $ or $ (\mathsf{B}) $. An example for the \textit{figure-eight knot} is shown in \Figs{Fig:figure8}.

\begin{figure}[!ht]
\centering
\hspace*{\fill}
\subfigure[The states of the figure-eight knot  following  the initial ``$ \mathsf{A}$'' split.]{\includegraphics[width=.925\linewidth]{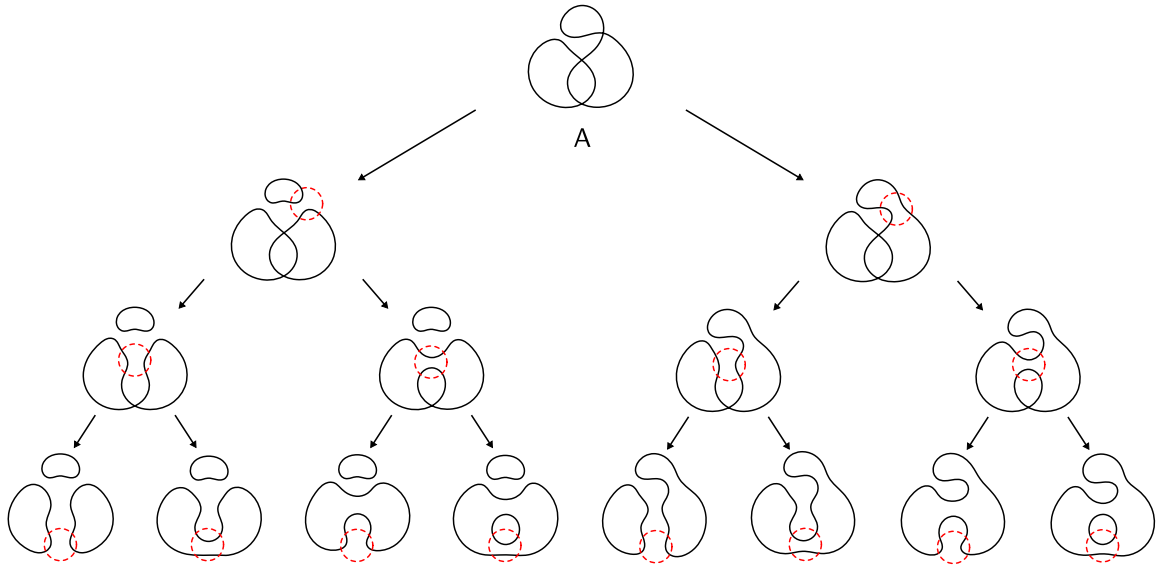}}
\hspace*{\fill}
\\
\hspace*{\fill}
\subfigure[The states of the figure-eight knot following  the initial  ``$ \mathsf{B}$'' split.]{\includegraphics[width=.925\linewidth]{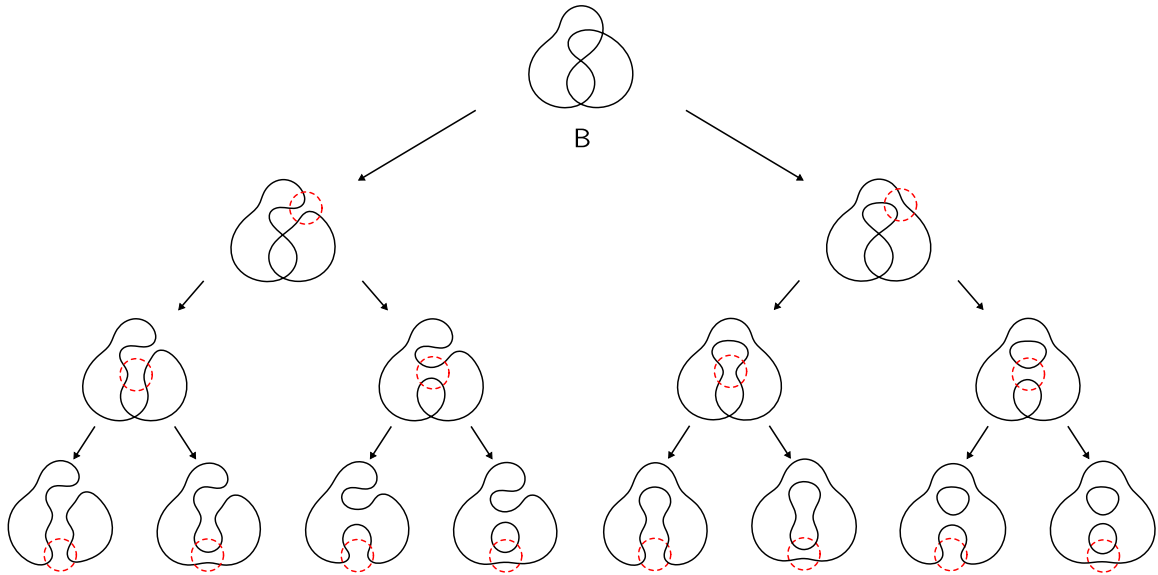}}
\hspace*{\fill}
\caption{The states of the figure-eight knot.}
\label{Fig:figure8states}
\end{figure}

Let $ K $ be a knot diagram. If $ m $ denotes the initial number of crossings, then the final states form a collection of $ 2^m $ diagrams of nonintersecting curves. We can enumerate those states with respect to the number of their components -- called \textit{circles} -- by introducing the sum
\begin{equation}\label{Eq:gp}
K(x):= \sum_{S}^{}x^{|S|},
\end{equation} where $ S $ browses the collection of  states, and $ |S| $ gives the number of circles in $ S $. Here, $ K(x) $ is an integer polynomial which we referred to as generating polynomial \cite{Ramaharo1,Ramaharo2} (in fact, it is a simplified formulation of what Kauffman calls  ``state polynomial'' \cite[Section 1--2]{Kauffman1} or    ``bracket polynomial'' \cite{Kauffman2}). For instance, if $ K $ is the figure-eight knot diagram, then we have $ K(x)=5x+8x^2+3x^3 $ (the states are illustrated in \Figs{Fig:figure8states}).

\begin{figure}[!ht]
\centering
\hspace*{\fill}
\subfigure[$B_{n,r}\rightleftharpoons B_{r,n}$.]{\includegraphics[width=.8\linewidth]{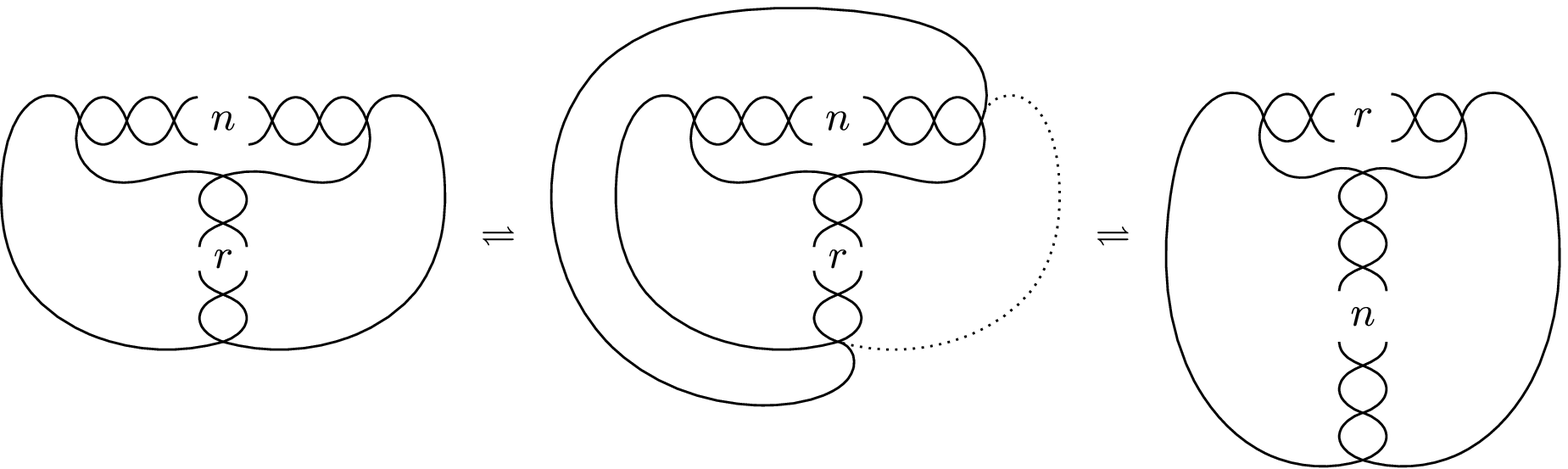}\label{subfig:2bnr}}
\hspace*{\fill}
\\
\subfigure[$B_{n,0} $.]{\includegraphics[width=0.255\linewidth]{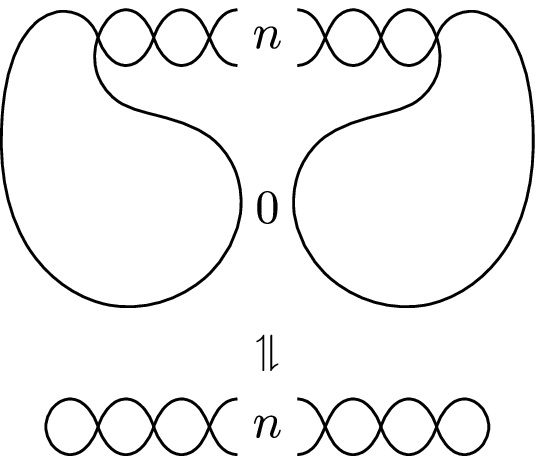}\label{subfig:2bn0}}\hfill%
\subfigure[$B_{n,\infty} $ \big($(2,n)  $-torus knot\big).]{\includegraphics[width=0.625\linewidth]{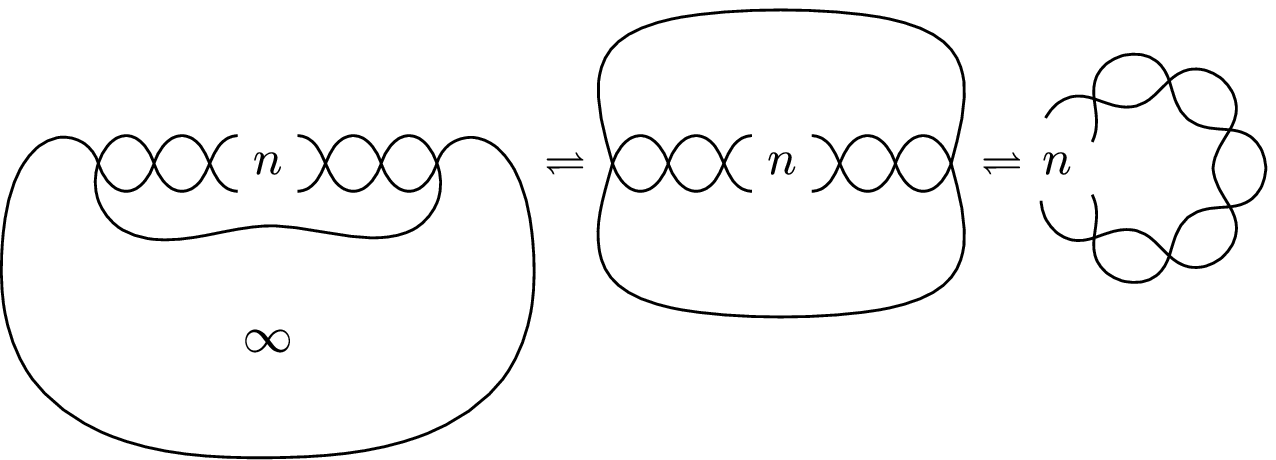}\label{subfig:2bninfty}}
\\
\subfigure[$B_{0,r}$.]{\includegraphics[width=0.325\linewidth]{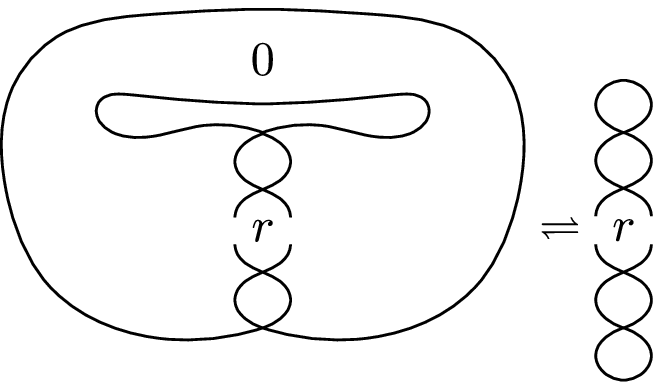}\label{subfig:2b0r}}\hfill %
\subfigure[$B_{\infty,r}$  \big($(2,r)  $-torus knot\big).]{\includegraphics[width=.55\linewidth]{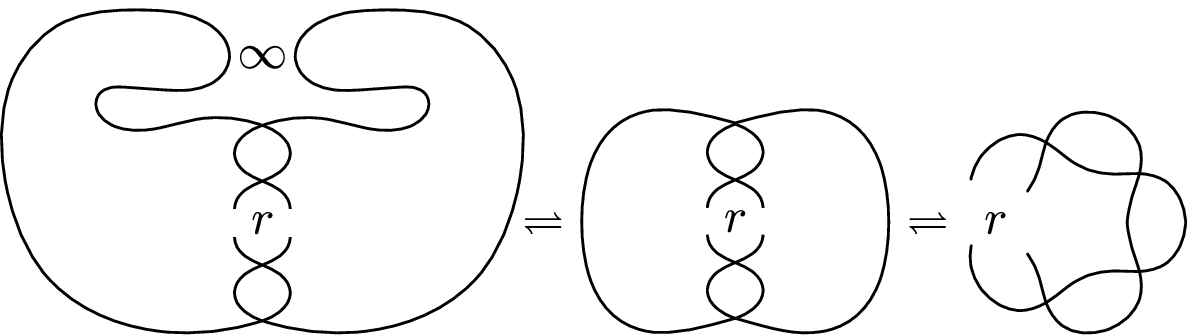}\label{subfig:2binftyr}}
\caption{Two-bridge knots with Conway's notation $ \mathcal{C}(n,r) $.}
\label{Fig:2bknot}
\end{figure}

In this note, we  establish the generating polynomial for the \textit{two-bridge knot}  with Conway's notation $ \mathcal{C}(n,r) $ \cite{Lafferty,Overduin}. We refer to the associated knot diagram as $ B_{n,r} $, where $ n $ and $ r $ denote the number of half-twists. For example, the figure-eight knot has Conway's notation $ \mathcal{C}(2,2) $.  Owing to the property of the shadow diagram which we draw on the sphere \cite{Denton}, we can continuously deform the diagram $ B_{n,r} $ into $ B_{r,n} $ without altering the crossings configuration. We let  $ B_{n,r} \rightleftharpoons B_{r,n} $ express such transformation \big(see \Figs{Fig:2bknot}~\subref{subfig:2bnr}\big). Besides, we let  $ B_{n,0} $ and $ B_{n,\infty} $ denote the diagrams in \Figs{Fig:2bknot}~\subref{subfig:2bn0} and \subref{subfig:2bninfty}, respectively. Here, ``$ 0 $'' and  ``$ \infty $'' are symbolic notations -- borrowed from tangle theory \cite[p.\ 88]{Kauffman1} -- that express the absence of half-twists. If $ r=\infty $ and $ n\geq 1 $, then $ B_{n,\infty}  $ represents the diagram of a $(2,n)  $-torus knot \big($ \rightleftharpoons B_{n-1,1} $\big).  Correspondingly, we let $ B_{0,r} $ and $ B_{\infty,r} $ denote the diagrams pictured in \Figs{Fig:2bknot}~\subref{subfig:2b0r} and \subref{subfig:2binftyr}, respectively.

\section{Generating polynomial}
Let $ K $, $ K' $ and $ \bigcirc $ be knot diagrams, where $ \bigcirc $ is the trivial knot, and let $ \# $ and $ \sqcup $ denote the connected sum and the disjoint union, respectively. The generating polynomial defined in \eqref{Eq:gp} verifies the following basic properties:
\begin{enumerate}[(i)]
\item $ \bigcirc(x)=x $\label{it:p1};	
\item $ \big(K\sqcup K'\big)(x) = K(x)K'(x) $\label{it:p2};
\item $ \big(K\# K'\big)(x) = \dfrac{1}{x}K(x)K'(x) $\label{it:p3}.
\end{enumerate}
Furthermore, if $ K\rightleftharpoons K' $, then $ K(x)=K'(x) $ \cite{Ramaharo1}.
\begin{lemma}\label{Lem:gpft}
The generating polynomial for the knots $ B_{n,0} $  and $ B_{n,\infty} $  are given by
\begin{align}\label{Eq:b0n}
B_{n,0}(x)=x(x+1)^n
\end{align}
and 
\begin{equation}\label{Eq:binftyn}
B_{n,\infty}(x)=(x+1)^{n} + x^2 -1.
\end{equation}
\end{lemma}
The key ingredient for establishing \eqref{Eq:b0n} and \eqref{Eq:binftyn} consists of the states of specific crossings which produce the recurrences 
\begin{equation*}
B_{n,0}(x)=\left(\bigcirc\sqcup B_{n-1,0}\right)(x)+B_{n-1,0}(x)
\end{equation*}
and
\begin{equation*}
B_{n,\infty}(x)=B_{n-1,0}(x)+B_{n-1,\infty}(x),
\end{equation*}
respectively, with initial values $ B_{0,0}(x)=x $ and $ B_{0,\infty}(x)=x^2 $ \cite{Ramaharo1}. Note that the lemma still holds if we replace index $ n $ by $ r $.

\begin{proposition}\label{Lem:closedform}
The generating polynomial for the two-bridge knot $ B_{n,r} $ is  given by the recurrence
\begin{equation}
B_{n,r}(x)=B_{n-1,r}(x)+(x+1)^{n-1}B_{\infty,r}(x),
\end{equation}
and has the following closed form:
\begin{equation}\label{Eq:closedform}
B_{n,r}(x)= \left(\dfrac{(x+1)^r+x^2-1}{x}\right)(x+1)^n+\left(x^2-1\right)\left(\dfrac{(x+1)^r-1}{x}\right).
\end{equation}

\begin{proof}
\begin{figure}[!ht]
\centering
\includegraphics[width=.95\linewidth]{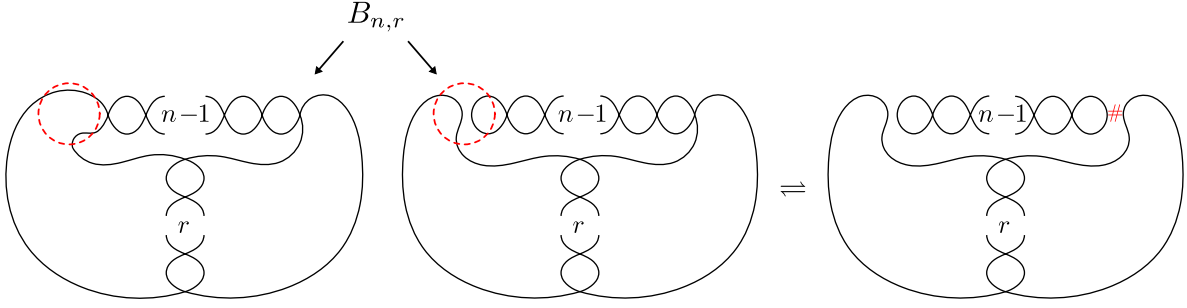}
\caption{The splits at a crossing allow us to capture $ B_{n-1,r} $, $ B_{n-1,0} $ and $ B_{\infty,r} $.}
\label{Fig:recBnr}
\end{figure}

By \Figs{Fig:recBnr} we have	
\begin{align*}
B_{n,r}(x)&=B_{n-1,r}(x)+\big(B_{n-1,0}\#B_{\infty,r}\big)(x)\\
&=B_{n-1,r}(x)+(x+1)^{n-1}B_{\infty,r}(x),
\end{align*}
where the last relation follows from property \eqref{it:p3}. Solving the recurrence for $ n $ yields
\[
B_{n,r}(x)=B_{0,r}(x)+B_{\infty,r}(x)\left(\dfrac{(x+1)^n-1}{x}\right).
\]
We conclude by the closed forms in \Lem{Lem:closedform}.
\end{proof}
\end{proposition}

\begin{remark}
We can write  \begin{equation}\label{Eq:bn0alpha}
B_{n,0}(x)=x^2\upalpha_n(x)+x 
\end{equation}
and
\begin{equation}\label{Eq:bninftyalpha}
B_{n,\infty}(x)=x\upalpha_n(x)+x^2,
\end{equation}
 where $ \upalpha_n(x):=\dfrac{(x+1)^n-1}{x} $, so that identity \eqref{Eq:closedform} becomes 
  \begin{equation}\label{Eq:Bnrx2x}
 B_{n,r}(x)=\left(x^2\upalpha_n(x) + x\right)+\left(x^2\upalpha_r(x)+x\upalpha_n(x)\upalpha_r(x)\right).
 \end{equation}
 
Since the coefficients of  $ \upalpha_n(x) $ are all nonnegative, it is clear, by \eqref{Eq:bn0alpha}, that the polynomial $ x^2\upalpha_n(x) $ counts the states of $ B_{n,0} $ that have at least $ 2 $ circles. This is illustrated in \Figs{Fig:statesbn0infty}~\subref{subfig:bn0}. Likewise, we have an interpretation of \eqref{Eq:bninftyalpha} in \Figs{Fig:statesbn0infty}~\subref{subfig:bninfty}. In \Figs{Fig:statesbn0infty} and \ref{Fig:partition}, the dashed diagrams represent all possible disjoint union of $ \ell-1 $ circles ($ \ell=n $ or $ r $, depending on the context), counted by $ \upalpha_\ell(x) $ and eventually empty. 
 
\begin{figure}[!ht]
\centering
\hspace*{\fill}
\subfigure[Illustration of $ x^2\upalpha_n(x)+x $.]{\includegraphics[width=.475\linewidth]{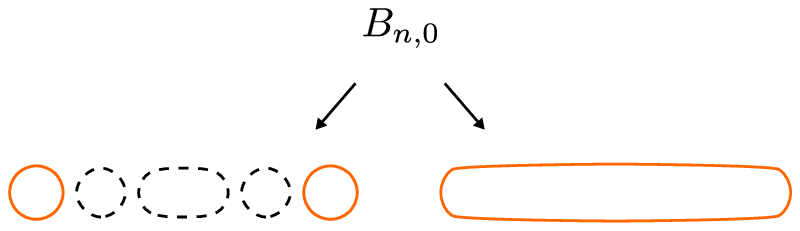}\label{subfig:bn0}}\hfill %
\subfigure[Illustration of $ x\upalpha_n(x)+x^2 $.]{\includegraphics[width=0.475\linewidth]{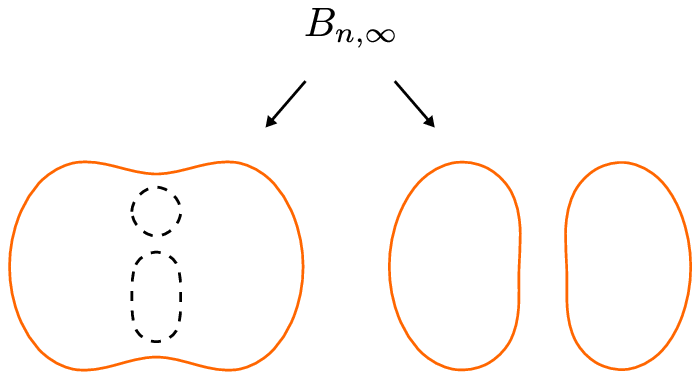}\label{subfig:bninfty}}\hfill %
\hspace*{\fill}
\caption{Illustrations of $ B_{n,0}(x) $ and $ B_{n,\infty}(x) $ as functions of $ \upalpha_n(x) $ .}
\label{Fig:statesbn0infty}
\end{figure}
\end{remark}

Therefore, for $ n,r\notin \left\{0,\infty\right\} $, identity \eqref{Eq:Bnrx2x}  means that we can classify the states into $ 4 $ subset as shown in \Figs{Fig:partition}. In these illustrations, there are $ 2^n-1 $ and $ 2^r-1 $  states of  \subref{subfig:s1}  and \subref{subfig:s2} kind, respectively, and $ \binom{n}{1}\times\binom{r}{1} +1$ one-component states of \subref{subfig:s3} and \subref{subfig:s4} kind. The remaining states are of  \subref{subfig:s3} kind, bringing the total number of states to $ 2^{n+r} $.

\begin{figure}[!ht]
\centering
\hspace*{\fill}
\subfigure[]{\includegraphics[width=.22\linewidth]{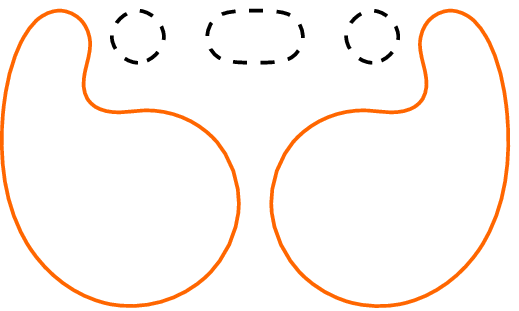}\label{subfig:s1}}\hfill %
\subfigure[]{\includegraphics[width=.22\linewidth]{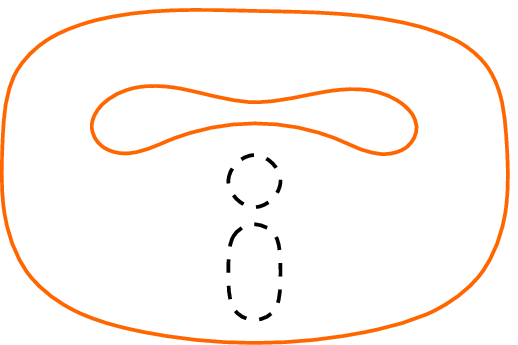}\label{subfig:s2}}\hfill %
\subfigure[]{\includegraphics[width=.22\linewidth]{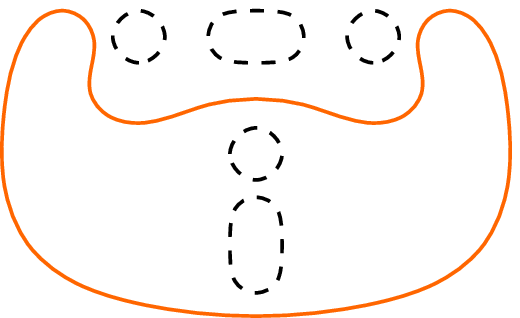}\label{subfig:s3}}\hfill %
\subfigure[]{\includegraphics[width=.22\linewidth]{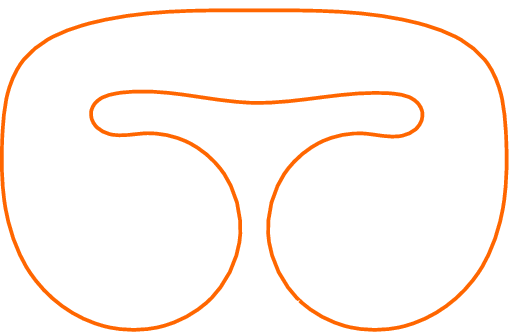}\label{subfig:s4}}\hfill %
\hspace*{\fill}
\caption{The states of $ B_{n,r} $: states in \subref{subfig:s1} are counted by $ x^2\upalpha_n(x) $, those in \subref{subfig:s2} by $ x^2\upalpha_r(x) $, those in \subref{subfig:s3} by $ x\upalpha_n(x)\upalpha_r(x) $, and state in \subref{subfig:s4} is simply counted by $ x $.}
\label{Fig:partition}
\end{figure}

\section{Particular values}
Let $\sum_{k\geq 0} b{(n,r;k)}x^k:=B_{n,r}(x)$, or  $ b(n,r;k):=\left[x^k\right] B_{n,r}(x) $. Then
\begin{equation*}
b(n,r;k)=\binom{n+r}{k+1}+\binom{n}{k-1}+\binom{r}{k-1}-\binom{n}{k+1}-\binom{r}{k+1}-\delta_{1,k},
\end{equation*}
where $ \delta_{1,k} $ is  the Kronecker symbol. By \eqref{Eq:gp}, we  recognize $ b(n,r;k) $ as the cardinal of the set $ \{|S|=k: S\textit{ is a state of }B_{n,r}\} $, i.e., the number of states having $ k $ circles. In this section, the coefficients $ b(n,r;k) $ are tabulated for some values of $ n $, $ r $ and $ k $. We give as well the corresponding A-numbers in the \textit{On-Line Encyclopedia of Integer Sequences} \cite{Sloane}.

\begin{itemize}
\item $ b(n,0;k)=\left[x^k\right] x(x+1)^n $, essentially giving entries in Pascal's triangle \seqnum{A007318}  (see \Tabs{Tab:bn0k}).
\begin{table}[!ht]
\centering
$\begin{array}{c|rrrrrrrrrrrrrrrrrrrrrrrrr}
n\ \backslash\ k		 &0		 &1		 &2		 &3		 &4		 &5		 &6		 &7	&8\\
\midrule
0 & 0 & 1\\
1 & 0 & 1 & 1\\
2 & 0 & 1 & 2 &  1\\
3 & 0 & 1 & 3 &  3 &  1\\
4 & 0 & 1 & 4 &  6 &  4 &  1\\
5 & 0 & 1 & 5 & 10 & 10 &  5 &  1\\
6 & 0 & 1 & 6 & 15 & 20 & 15 &  6 & 1\\
7 & 0 & 1 & 7 & 21 & 35 & 35 & 21 & 7 & 1
\end{array}$
\caption{Values of $ b(n,0;k) $ for $ 0\leq n\leq 7 $ and $ 0\leq k\leq 8 $.}
\label{Tab:bn0k}
\end{table}
	
\item $ b(n,1;k)=\left[x^k\right]\left((x+1)^{n+1}+x^2-1\right) $, generating a subtriangle in \seqnum{A300453} (see \Tabs{Tab:bn1k}).
\begin{table}[!ht]
\centering
$\begin{array}{c|rrrrrrrrrrrrrrrrrrrrrrrrr}
n\ \backslash\ k		 &0		 &1		 &2		 &3		 &4		 &5		 &6		 &7	&8\\
\midrule
0 & 0 & 1 &  1\\
1 & 0 & 2 &  2\\
2 & 0 & 3 &  4 &  1\\
3 & 0 & 4 &  7 &  4 &  1\\
4 & 0 & 5 & 11 & 10 &  5 &  1\\
5 & 0 & 6 & 16 & 20 & 15 &  6 &  1\\
6 & 0 & 7 & 22 & 35 & 35 & 21 &  7 & 1\\
7 & 0 & 8 & 29 & 56 & 70 & 56 & 28 & 8 & 1\\
\end{array}$
\caption{Values of $ b(n,1;k) $ for $ 0\leq n\leq 7 $ and $ 0\leq k\leq 8 $.}
\label{Tab:bn1k}
\end{table}
	
\item $ b(n,2;k)=\left[x^k\right]\left((2x+2)(x+1)^{n}+\left(x^2-1\right)(x+2)\right) $, giving triangle in \seqnum{A300454} (see \Tabs{Tab:bn2k}).
\begin{table}[!ht]
\centering
$\begin{array}{c|rrrrrrrrrrrrrrrrrrrrrrrrr}
n\ \backslash\ k		 &0		 &1		 &2		 &3		 &4		 &5		 &6		 &7		 &8\\
\midrule
0 & 0 &  1 &  2 &   1\\
1 & 0 &  3 &  4 &   1\\
2 & 0 &  5 &  8 &   3\\
3 & 0 &  7 & 14 &   9 &   2\\
4 & 0 &  9 & 22 &  21 &  10 &   2\\
5 & 0 & 11 & 32 &  41 &  30 &  12 &  2\\
6 & 0 & 13 & 44 &  71 &  70 &  42 & 14 &  2\\
7 & 0 & 15 & 58 & 113 & 140 & 112 & 56 & 16 & 2
\end{array}$
\caption{Values of $ b(n,2;k) $ for $ 0\leq n\leq 7 $ and $ 0\leq k\leq 8 $.}
\label{Tab:bn2k}
\end{table}

\item $ b(n,n;k)=\left[x^k\right] \left(\dfrac{(x + 1)^{2n} + \big(x^2 - 1\big)\big(2(x + 1)^n - 1)\big)}{x}\right)$, giving triangle in \seqnum{A321127} (see \Tabs{Tab:bnnk}).
\begin{table}[!ht]
\centering
$\begin{array}{c|rrrrrrrrrrrrrrrrrrrrrrrrr}
n\ \backslash\ k		 &0		 &1		 &2		 &3		 &4		 &5		 &6		 &7		 &8		 &9		 &10		 &11		 &12		 &13\\
\midrule
0 & 0 &  1\\
1 & 0 &  2 &   2\\
2 & 0 &  5 &   8 &   3\\
3 & 0 & 10 &  24 &  21 &    8 &    1\\
4 & 0 & 17 &  56 &  80 &   64 &   30 &    8 &    1\\
5 & 0 & 26 & 110 & 220 &  270 &  220 &  122 &   45 &   10 &    1\\
6 & 0 & 37 & 192 & 495 &  820 &  952 &  804 &  497 &  220 &   66 &  12 &  1\\
7 & 0 & 50 & 308 & 973 & 2030 & 3059 & 3472 & 3017 & 2004 & 1001 & 364 & 91 & 14 & 1
\end{array}$
\caption{Values of $ b(n,n;k) $ for $ 0\leq n\leq 7 $ and $ 0\leq k\leq 13 $.}
\label{Tab:bnnk}
\end{table}

\item $ b(n,r;1) = nr + 1  $, giving \seqnum{A077028}, and displayed as square array in \Tabs{Tab:bnr1}.
\begin{table}[!ht]
\centering
$\begin{array}{c|rrrrrrrrrrrrrrrrrrrrrrrrr}
n\ \backslash\ r		 &0		 &1		 &2		 &3		 &4		 &5		 &6		 &7\\
\midrule
0 & 1 &  1 &  1 &  1 &  1 &  1 &  1 &  1\\
1 & 1 &  2 &  3 &  4 &  5 &  6 &  7 &  8\\
2 & 1 &  3 &  5 &  7 &  9 & 11 & 13 & 15\\
3 & 1 &  4 &  7 & 10 & 13 & 16 & 19 & 22\\
4 & 1 &  5 &  9 & 13 & 17 & 21 & 25 & 29\\
5 & 1 &  6 & 11 & 16 & 21 & 26 & 31 & 36\\
6 & 1 &  7 & 13 & 19 & 25 & 31 & 37 & 43\\
7 & 1 &  8 & 15 & 22 & 29 & 36 & 43 & 50
\end{array}$
\caption{Values of $ b(n,r;1) $ for $ 0\leq n\leq 7 $ and $ 0\leq r\leq 7 $.}
\label{Tab:bnr1}
\end{table}

In Kauffman's language, $ b(n,r;1) $ is, for a fixed choice of star region, the number of ways of placing \textit{state markers} at the crossings of the diagram $ B_{n,r} $, i.e., of the forms 
\[
\protect\includegraphics[width=0.06\linewidth,valign=c]{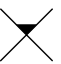},\quad
\protect\includegraphics[angle=90,origin=c,width=0.06\linewidth,valign=c]{markersA1},\quad
\protect\includegraphics[angle=180,origin=c,width=0.06\linewidth,valign=c]{markersA1},\quad
\protect\includegraphics[angle=270,origin=c,width=0.06\linewidth,valign=c]{markersA1},
\] so that the resulting states are  ``Jordan trails'' \cite[Section 1--2]{Kauffman1}. Note that a state marker is interpreted as an instruction to split a crossing as shown below: 
\[
\protect\includegraphics[width=0.06\linewidth,valign=c]{markersA1}\Longrightarrow\protect\includegraphics[width=0.06\linewidth,valign=c]{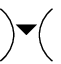}\qquad\mbox{ and }\qquad\protect\includegraphics[angle=90,origin=c,width=0.06\linewidth,valign=c]{markersA1}\Longrightarrow\protect\includegraphics[angle=90,origin=c,width=0.06\linewidth,valign=c]{markersA2}.
\]
The process is illustrated in \Figs{Fig:statemarkers} for the figure-eight knot. 
\begin{figure}[!ht]
\centering
\includegraphics[width=.7\linewidth]{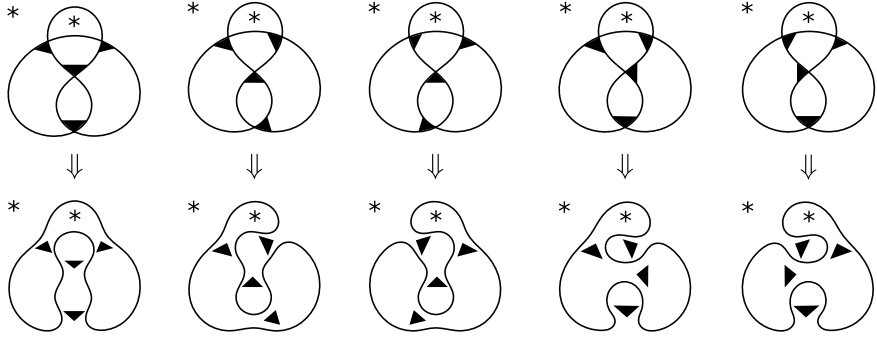}
\caption{Illustration of $ b(2,2;1) $: mark two adjacent regions by stars (\textasteriskcentered), then assign a state marker at each crossing so that no region of $B_{2,2}$ contains more than one state marker, and regions with stars do not have any.}
\label{Fig:statemarkers}
\end{figure}

\item $ b(n,r;2) = n\left(\binom{r}{2}+1\right)+r\left(\binom{n}{2}+1\right) $, giving square array in \seqnum{A300401} (see \Tabs{Tab:bnr2}).
\begin{table}[!ht]
\centering
$\begin{array}{c|rrrrrrrrrrrrrrrrrrrrrrrrr}
n\ \backslash\ r		 &0		 &1		 &2		 &3		 &4		 &5		 &6		 &7\\
\midrule
0 &  0 &  1 &   2 &   3 &   4 &   5 &   6 &   7\\
1 &  1 &  2 &   4 &   7 &  11 &  16 &  22 &  29\\
2 &  2 &  4 &   8 &  14 &  22 &  32 &  44 &  58\\
3 &  3 &  7 &  14 &  24 &  37 &  53 &  72 &  94\\
4 &  4 & 11 &  22 &  37 &  56 &  79 & 106 & 137\\
5 &  5 & 16 &  32 &  53 &  79 & 110 & 146 & 187\\
6 &  6 & 22 &  44 &  72 & 106 & 146 & 192 & 244\\
7 &  7 & 29 &  58 &  94 & 137 & 187 & 244 & 308\\
\end{array}$
\caption{Values of $ b(n,r;2) $ for $ 0\leq n\leq 7 $ and $ 0\leq r\leq 7 $.}
\label{Tab:bnr2}
\end{table}

We paid a special attention to the case $ k=2 $ because, surprisingly, columns $\big(b(n,1;2)\big)_n$ and $\big(b(n,2;2)\big)_n$ match sequences \seqnum{A000124} and \seqnum{A014206}, respectively \cite{Ramaharo1}. The former gives the maximal number of regions into which the plane is divided by $ n $ lines, and the latter the maximal number of regions into which the plane is divided by $ (n+1) $ circles.

\item $ b\big(n,r;d(n,r)\big)=\textit{leading coefficient of }B_{n,r}(x)$, giving square array in \seqnum{A321125} (see \Tabs{Tab:bnrdnr}). Here,  $ d(n,r)=\max(n+1,r+1,n+r-1) $ denotes the degree of $ B_{n,r}(x) $, and gives entries in \seqnum{A321126}. We have \Tabs{Tab:degree} giving the numbers $ d(n,r) $ for $ 0\leq n\leq 7 $ and $ 0\leq r\leq 7 $. 
\begin{table}[!ht]
\centering
$\begin{array}{c|rrrrrrrrrrrrrrrrrrrrrrrrr}
n\ \backslash\ r		 &0		 &1		 &2		 &3		 &4		 &5		 &6 &7\\
\midrule
0 & 1 &  1 &  1 &  1 &  1 &  1 &  1 &  1\\
1 & 1 &  2 &  1 &  1 &  1 &  1 &  1 &  1\\
2 & 1 &  1 &  3 &  2 &  2 &  2 &  2 &  1\\
3 & 1 &  1 &  2 &  1 &  1 &  1 &  1 &  1\\
4 & 1 &  1 &  2 &  1 &  1 &  1 &  1 &  1\\
5 & 1 &  1 &  2 &  1 &  1 &  1 &  1 &  1\\
6 & 1 &  1 &  2 &  1 &  1 &  1 &  1 &  1\\
7 & 1 &  1 &  2 &  1 &  1 &  1 &  1 &  1
\end{array}$
\caption{Leading coefficients of $ B_{n,r}(x) $ for $ 0\leq n\leq 7 $ and $ 0\leq r\leq 7 $.}
\label{Tab:bnrdnr}
\end{table}

\begin{table}[!ht]
\centering
$\begin{array}{c|rrrrrrrrrrrrrrrrrrrrrrrrr}
n\ \backslash\ r		 &0		 &1		 &2		 &3		 &4		 &5		 &6		 &7\\
\midrule
0 & 1 & 2 & 3 & 4 &  5 &  6 &  7 &  8\\
1 & 2 & 2 & 3 & 4 &  5 &  6 &  7 &  8\\
2 & 3 & 3 & 3 & 4 &  5 &  6 &  7 &  8\\
3 & 4 & 4 & 4 & 5 &  6 &  7 &  8 &  9\\
4 & 5 & 5 & 5 & 6 &  7 &  8 &  9 & 10\\
5 & 6 & 6 & 6 & 7 &  8 &  9 & 10 & 11\\
6 & 7 & 7 & 7 & 8 &  9 & 10 & 11 & 12\\
7 & 8 & 8 & 8 & 9 & 10 & 11 & 12 & 13\\
\end{array}$
\caption{Values of $ d(n,r) $ for $ 0\leq n\leq 7 $ and $ 0\leq r\leq 7 $.}
\label{Tab:degree}
\end{table}

We have the following properties:
\begin{itemize}
\item $ d(n,r)=d(r,n)$;
\item if $ r=0 $, then $ d(n,r)=n+1 $;
\item if $ r=\infty $, then sequence $ \big(d(n,r)\big)_n$ begins: $ 2,2,2,3,4,5,6,7,8,\ldots $ (\seqnum{A233583} with offset 1).
\end{itemize}

Diagramatically, we give the corresponding illustration  for some values of $ n $ and $ r $  in \Figs{Fig:dnr}.
\begin{figure}[!ht]
\centering
\hspace*{\fill}
\subfigure[$ n\geq 1 $.]{\includegraphics[width=.25\linewidth]{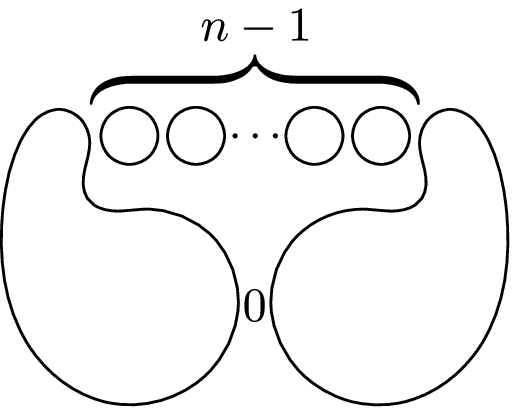}\label{subfig:dn0}}\hfill %
\subfigure[$ n\geq 2 $.]{\includegraphics[width=.25\linewidth]{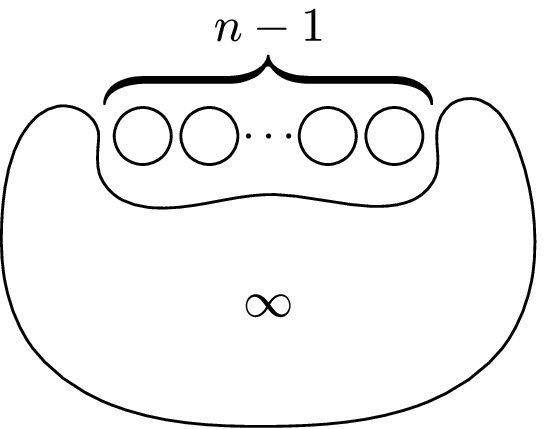}\label{subfig:dninfty}}\hfill %
\subfigure[$ n\geq 2 $, $ r\geq 2 $.]{\includegraphics[width=0.25\linewidth]{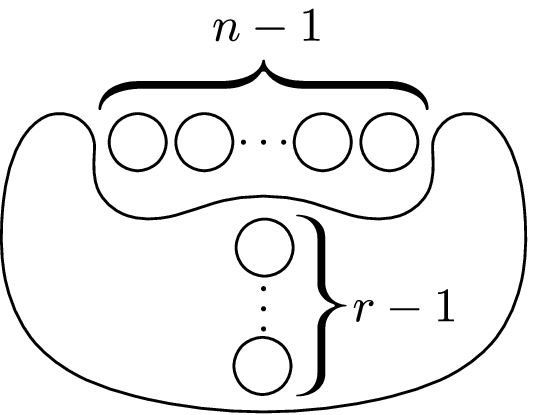}\label{subfig:dnr}}\hfill %
\hspace*{\fill}
\caption{Illustration of the numbers $ d(n,r) $.}
\label{Fig:dnr}
\end{figure}

Correspondingly, we have
\begin{itemize}
\item $ b\big(n,r;d(n,r)\big)=b\big(r,n;d(r,n)\big)$;
\item if $ r=0 $, then $ b\big(n,r;d(n,r)\big)=1 $;
\item if $ r=\infty $, then sequence $ \Big(b\big(n,r;d(n,r)\big)\Big)_n$ begins: $ 1,1,2,1,1,1,1,1,\ldots $ (\seqnum{A294619} with initial term equals to $ 0 $).
\end{itemize}

Remarquable values in \Tabs{Tab:bnrdnr} correspond to knots $ B_{1,1} $ (``Hopf link'', see \Figs{Fig:hopfstates}), $ B_{2,2} $ (figure-eight knot, see \Figs{Fig:figure8}, \ref{Fig:figure8states}) and $ B_{n,2} $ (``twist knot'' \cite{Ramaharo1}) for $ n\geq 3 $. The latter case can be observed from identity \eqref{Eq:Bnrx2x} for which the leading coefficient is larger than $ 1 $ when $ n+1=n+r-1 $ is satisfied. Also, considere the identity below:
\[
B_{n,2}(x) =  B_{n,0}(x) + B_{n,\infty}(x) + B_{n,\infty}(x) + \big(\bigcirc\sqcup B_{n,\infty}\big)(x).
\]
We notice that the leading coefficient is deduced from the contribution of the summands $ B_{n,0}(x) $ and $ \big(\bigcirc\sqcup B_{n,\infty}\big)(x) $ \cite{Ramaharo1}.
\begin{figure}[!ht]
\centering
\includegraphics[width=.5\linewidth]{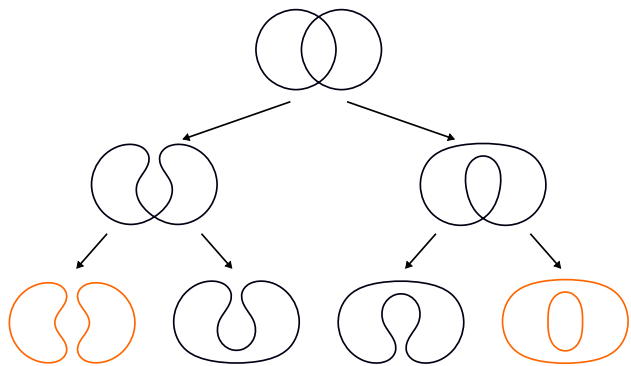}
\caption{The states of the knot  $ B_{1,1} $: $ d(1,1)=2 $ and $ b\big(1,1,d(1,1)\big)=2 $.}
\label{Fig:hopfstates}
\end{figure}
\end{itemize}

\bigskip
\hrule
\bigskip

\noindent 2010 {\it Mathematics Subject Classification}: 05A10; 57M25.

\bigskip
\hrule
\bigskip

\noindent (Concerned with sequences
\seqnum{A000124}, \seqnum{A007318}, \seqnum{A014206}, \seqnum{A077028}, \seqnum{A233583}, \seqnum{A294619},\\ \seqnum{A300401}, \seqnum{A300453}, \seqnum{A300454}, \seqnum{A321125}, \seqnum{A321126} and \seqnum{A321127}.)
\bigskip
\hrule
\bigskip
		
\end{document}